\documentclass{amsart}
\usepackage{amssymb}
\usepackage{amscd}
\usepackage{graphics}

\usepackage[
hyperref,
nonotes
]{degt}

\def\CG#1{\Z/#1}

\let\Fermat=\Phi
\let\PQ=\Sigma
\let\rl=R

\let\ln=L
\def\tln{\reduced\ln}
\let\dv=V

\def\sset{\bold S}
\def\ofD#1{\<#1\>}
\def\SS{\sset\ofD}
\def\KK{\bold K\ofD}
\def\TT{\bold T\ofD}
\def\QQ{\bold Q}
\def\GG{\Bbb G}
\def\Pic{\QOPNAME{Pic}}
\def\Tor{\QOPNAME{Tor}}
\def\CO{\Cal O}

\def\of#1{(#1)}
\def\of#1{[#1]}
\def\ofm{\of{m}}

\def\tS{\tilde S}
\def\tU{\tilde U}

\let\Am=A
\let\Hm=A
\def\Bm{\Am^\circ}
\def\tAm{\reduced\Am}

\let\=\B

\let\vg=v
\let\hg=h
\let\ag=a
\let\cg=c
\let\ug=u
\let\g=g
\let\3\tilde
\def\bag{\3\ag}
\def\bbg{\3\cg}
\let\GL\Lambda
\let\reduced\bar
\def\tGL{\reduced\GL}

\def\longinto{\DOTSB\lhook\joinrel\longrightarrow}

\def\bGa{\tilde\Ga}

\let\poly\varphi
\let\inc\iota

\def\CW{{\sl CW}}

\title[Lines generate the Picard groups]
 {Lines generate the Picard groups\\of certain Fermat surfaces}

\author{Alex Degtyarev}

\address{%
Department of Mathematics\\
Bilkent University\\
06800 Ankara, Turkey}

\email{degt@fen.bilkent.edu.tr}



\begin{document}

\begin{abstract}
We answer a question of T.~Shioda and show that, for any positive integer~$m$
prime to~$6$, the Picard group of the Fermat surface $\Fermat_m$ is generated
by the classes of lines contained in~$\Fermat_m$.
A few other classes of surfaces are also considered.
\end{abstract}

\maketitle

\section{Introduction}

\subsection{Principal results}
All algebraic varieties in the paper are over~$\C$.
Let~$m$ be a positive integer, and let
\[*
\Fermat_m:=\{z_0^m+z_1^m+z_2^m+z_3^m=0\}\subset\Cp3
\]
be the \emph{Fermat surface}.
If $m=1$ (plane) or $m=2$ (quadric), then\mnote{`then'}
$\Fermat_m$ contains
infinitely many lines (meaning true straight lines in~$\Cp3$);
otherwise, $\Fermat_m$ is known to contain exactly $3m^2$ lines.

Since $\Fermat_m$ is simply connected, one can identify its Picard group
$\Pic\Fermat_m$ and its N\'{e}ron--Severi lattice $\NS(\Fermat_m)$.
Citing~\cite{Shioda:Fermat1},
the N\'{e}ron--Severi group
``\dots is a rather delicate invariant of arithmetic nature.
Perhaps for this reason
it usually requires some nontrivial work before one can determine the Picard
number of a given variety, let alone the full structure of its N\'{e}ron--Severi
group.''
The Picard groups of Fermat surfaces are related to those of the more general
\emph{Delsarte surfaces}
(see~\cite{Shioda:Delsarte}; they fit into the framework outlined
in \autoref{s.m->1}).
Furthermore, continuing the citation,
``Combined with the method based on the inductive structure of Fermat
varieties, this might lead to the verification of the Hodge conjecture for
all Fermat varieties.''

Let $\sset_m\subset\Pic\Fermat_m$ be the subgroup generated by the classes of
the lines
contained
in~$\Fermat_m$. Then, according to~\cite{Shioda:Fermat}, one has
\[
\sset_m\otimes\Q=(\Pic\Fermat_m)\otimes\Q\quad
\text{if and only if $m\le4$ or $\gcd(m,6)=1$}.
\label{eq.Shioda}
\]
This statement is proved by comparing the dimensions of the two spaces, which
are computed independently.
In other words, the classes of lines generate $\Pic\Fermat_m$ rationally, and
a natural question, raised in~\cite{Shioda:Fermat1}, is whether they also
generate the Picard group over the integers.
A partial answer to this question
was given in~\cite{Shioda:Fermat2},
almost $30$
years later: the equality $\Pic\Fermat_m=\sset_m$ holds for
all integers~$m$ prime to~$6$ in
the range $5\le m\le100$.
This fact is proved by supersingular reduction and a computer aided
computation of the discriminants of the lattices involved.
(The case
$m=3$ is classical: any nonsingular cubic contains $27$ lines,
which generate
its Picard group.
The case $m=4$, \ie, that of $K3$-surfaces, was settled in~\cite{Muzukami},
see also~\cite{Boissiere.Sarti} for a slight generalization.
The proof suggested below works for both cases.)

The principal result of the present paper is the following theorem,
answering the
above question in the affirmative in the general case.

\theorem\label{th.Pic=S}
Let $m\ge1$ be an integer such that either $m\le4$ or $\gcd(m,6)=1$. Then
$\Pic\Fermat_m=\sset_m$, \ie,
$\Pic\Fermat_m$ is generated by the classes of lines.
\endtheorem

Since the $3m^2$ lines in~$\Fermat_m$ admit a very explicit description
(\cf. \autoref{s.m->1}) and one can easily see how they intersect
(see, \eg, Equation~(6) in~\cite{Shioda:Fermat2}),
\autoref{th.Pic=S} gives us a complete description of
$\Pic\Fermat_m=\NS(\Fermat_m)$,
including the intersection form and the action of the
automorphism group of~$\Fermat_m$.

In view of~\eqref{eq.Shioda},
\autoref{th.Pic=S} is
an immediate
consequence of the following statement, which is actually proved in the
paper, see \autoref{proof.main}.
(Throughout the paper, we use $\Tors A$ for the $\Z$-torsion of an abelian
group/module~$A$.)

\theorem\label{th.main}
For any integer $m\ge1$, one has $\Tors(\Pic\Fermat_m/\sset_m)=0$.
\endtheorem

In the mean time,
an interesting generalization,
approaching the problem from a different angle,
was suggested in~\cite{Shimada.Takahashi:primitivity}.
Briefly, $\Fermat_m$ can be represented as an $m^3$-fold
ramified covering of the plane, and one can try to study
other multiple planes with the same ramification locus (see \autoref{s.m->1}
and \autoref{problem} for details).
Considered in~\cite{Shimada.Takahashi:primitivity} are cyclic coverings of
degree at most~$50$, and, similar to~\cite{Shioda:Fermat2},
the proof is also based on comparing the
discriminants of the two lattices.

The
approach developed in the present paper,
including the computation of the Alexander
module~$\Hm\of\Ga$ (see \autoref{s.others}),
which is crucial for the proof, apply to Delsarte surfaces as well.
Here,
we make a few first steps towards this generalized problem and
work out another special case, see \autoref{th.toy}.
In the forthcoming paper~\cite{degt:Delsarte},
we close the case of \emph{cyclic} Delsarte surfaces started
in~\cite{Shimada.Takahashi:primitivity} and modify part of the proof of
\autoref{th.main} (see \autoref{proof.main}) to adapt it to slightly more
general \emph{diagonal} Delsarte surfaces.
On the other hand, numeric experiments show that \autoref{th.main} does not
extend literally to all Delsarte surfaces: sometimes, the quotient does have
torsion.
Next special classes to be studied
in more details would probably be the nonsingular
Delsarte
surfaces and those with
$\bA$--$\bD$--$\bE$ singularities.

As
yet another application,
we consider another class of surfaces whose Picard group is
rationally generated by lines, see~\cite{Boissiere.Sarti}.
Let $p$ and~$q$ be two square free
bivariate homogeneous polynomials of
degree~$m$, and denote
\[*
\PQ_{p,q}:=\{p(z_0,z_1)=q(z_2,z_3)\}\subset\Cp3.
\]
This nonsingular
surface contains an obvious set of $m^2$ lines, \viz. those connecting
the points $[z_0:z_1:0:0]$ and $[0:0:z_2:z_3]$, where
$p(z_0,z_1)=q(z_2,z_3)=0$, and we denote by
$\sset_{p,q}\subset\Pic\PQ_{p,q}$ the subgroup generated by
the classes of
these lines.

\theorem[see \autoref{proof.pq}]\label{th.pq}
For any pair $p$, $q$ as above, $\Tors(\Pic\PQ_{p,q}/\sset_{p,q})=0$.
\endtheorem

\corollary[see \autoref{proof.pq}]\label{cor.pq}
If $m$ is prime and $p$, $q$ as above are
sufficiently generic,
then $\Pic\PQ_{p,q}$ is
generated by the classes of the $m^2$ lines contained in $\PQ_{p,q}$.
\endcorollary

\subsection{An outline of the proof}\label{s.outline}
In \autoref{S.prelim}, we reduce the question to the computation of the
torsion of the $1$-homology of a certain space, see \autoref{th.reduction}.
We also recall the classical description of the lines in~$\Fermat_m$
by means of a ramified covering of the plane and,
following~\cite{Shimada.Takahashi:primitivity}, describe a generalization of
the problem to a wider class of surfaces.
In \autoref{S.Am}, we compute the so-called \emph{Alexander module} (or rather
Alexander complex) $\Hm\ofm$ of the above covering and its reduced version
$\tAm\ofm$.
The heart of the proof is a tedious computation of the length
$\ell(\tAm\ofm)$, see \autoref{lem.ell} in \autoref{S.proof};
then, \autoref{th.main} follows from
comparing the result to the expected value given
by~\cite{Shioda:Fermat1,Shioda:Fermat}, see \autoref{proof.main}.
In \autoref{s.toy}, we work out a toy example, illustrating the suggested
line of attack to the generalized problem.

\subsection{Acknowledgements}
I would like to express my gratitude to I.~Shimada for
bringing the problem to my attention and for many fruitful discussions;
it was he who eventually persuaded me to publish these observations.
I\mnote{added}
would also like to thank the anonymous referee of this paper for the elegant
proof of~\eqref{eq.ref}.

\section{Preliminaries}\label{S.prelim}

\subsection{Prerequisites}\label{s.alg.top}
For the reader's convenience,
we recall, with references to~\cite{Dold}, a few necessary facts
from algebraic topology. An ultimate reference would
be~\cite{Fomenko.Fuks};
unfortunately
it is only available in Russian.

By definition,
for any topological pair $(X,A)$ we have the following short
exact sequence of singular chain complexes:
\[*
0\longto S_*(A)\longto S_*(X)\longto S_*(X,A)\longto0.
\]
All complexes are free; hence, applying $\otimes\,G$ or $\Hom(\cdot,G)$, we
also have short exact sequences of (co-)chain complexes with any coefficient
group~$G$. The associated long exact sequences in (co-)homology are called
the (co-)homology exact sequences of pair $(X,A)$, \cf. $(3.2)$
in~\cite[Chapter III]{Dold}.

Unless specified otherwise, all (co-)homology
are with coefficients in~$\Z$.
The
other groups can be computed using the so-called
universal coefficient formulas
(see, \eg, $(7.9)$ and $(7.10)$ in \cite[Chapter VI]{Dold}):
for any topological space~$X$, abelian group~$G$, and integer~$n$, there are
natural split (not naturally) exact sequences
\[*
\gathered
0\longto H_n(X)\otimes G\longto
 H_n(X;G)\longto\Tor(H_{n-1}(X),G)\longto0,\\
0\longto\Ext(H_{n-1}(X),G)\longto
 H^n(X;G)\longto\Hom(H_n(X),G)\longto0.\\
\endgathered
\]
(Here, $\Tor=\Tor_1$ and $\Ext=\Ext^1$
are the derived functors in the category of $\Z$-modules.)
Similar statements hold for the relative groups of pairs $(X,A)$.
Assuming all groups finitely generated (\eg, $X$ is a finite \CW-complex),
a consequence of the second exact sequence is the assertion that
$H^n(X)$ is free if and only if so is $H_{n-1}(X)$; in this case,
$H^n(X)=\Hom(H_n(X),\Z)$.

We
use the following terminology for various duality isomorphisms
in topology of manifolds. Let $M$ be an oriented
compact manifold, $\dim M=n$, and
$A\subset M$ a
`sufficiently good'
(see the end of this paragraph)\mnote{edited}
closed subset. If $\partial M=\varnothing$, the multiplication
by the fundamental class $[M]$ establishes canonical isomorphisms
\roster*
\item
$H^p(M)=H_{n-p}(M)$ (Poincar\'{e} duality) and
\item
$H^p(M,A)=H_{n-p}(M\sminus A)$ (Poincar\'{e}--Lefschetz duality).
\endroster
In general, the multiplication by $[M,\partial M]$ is an isomorphism
\roster*
\item
$H^p(M)=H_{n-p}(M,\partial M)$ (Lefschetz duality).
\endroster
All statements are classical and well known.
For example, they can be derived as special
cases of Proposition~7.2 in \cite[Chapter VIII]{Dold},
with an extra observation that, in all cases considered in the paper, $M$ and
$A$ are at worst compact semialgebraic sets, thus admitting finite
triangulations (see, \eg, \cite{Hironaka.triang}); hence, they are absolute
neighborhood retracts and
the \v{C}ech cohomology in \cite{Dold} can be
replaced with singular ones.
As another consequence of~\cite{Hironaka.triang}, all (co-)homology groups
involved
are finitely generated.

\subsection{Divisors}\label{s.divisors}
Consider a
smooth projective
algebraic surface~$X$.
By Poincar\'{e} duality $H^2(X)=H_2(X)$, we can
regard the N\'{e}ron--Severi lattice $\NS(X)$ as a subgroup of
the intersection index lattice $H_2(X)/\Tors$,
representing a divisor $D\subset X$ by its fundamental class
$[D]$, see \autoref{s.homology} below.
(The N\'{e}ron--Severi \emph{lattice} is the group of divisors
modulo numeric equivalence; thus, we ignore the torsion.)
Since $\Pic X=H^1(X;\CO_X^*)$ and $H^2(X;\CO_X)$ is a $\C$-vector space,
the exponential exact sequence
\[
H^1(X;\CO_X)\longto H^1(X;\CO_X^*)\longto H^2(X)\longto H^2(X;\CO_X)
\label{eq.exp}
\]
implies that $\NS(X)$ is a primitive subgroup in $H_2(X)/\Tors$.

If $H_1(X)=0$, then
$H^2(X)=\Hom(H_2(X),\Z)$ is torsion free (the universal coefficient formula),
and so is $H_2(X)=H^2(X)$.
Since also $H^1(X;\CO_X)=H^{0,1}(X)$ is trivial in this case,
from~\eqref{eq.exp} we have $\Pic X=\NS(X)$, \ie, we do not need to
distinguish between linear, algebraic, or numeric equivalence of divisors.

Consider a reduced curve
$D\subset X$. Topologically, the normalization $\tilde D$ of~$D$ is a closed
surface,
and the projection $\Gs\:\tilde D\to D$ is a
homeomorphism outside a \emph{finite} subset $S\subset\tilde D$.
Hence,
\[*
H_2(D)=H_2(D,\Gs(S))=H_2(\tilde D,S)=H_2(\tilde D)=\bigoplus\Z\cdot[D_i]
\]
is the free abelian group generated by
the fundamental classes $[D_i]$ of
the irreducible components $D_i$ of~$D$ (or, equivalently,
the fundamental classes $[\tilde D_i]$ of the connected components
$\tilde D_i$ of~$\tilde D$). A similar computation in cohomology
shows that the group
\[*
H^2(D)=H^2(\tilde D)=\Hom(H_2(D),\Z)=\bigoplus\Z\cdot[D_i]^*
\]
is also free (the last identification uses the \emph{canonical} basis
$\{[D_i]\}$) and, by the universal coefficient formula, $H_1(D)$ is free.
(Essentially, we only use the fact that the singular locus has real
codimension at least two.)

\subsection{Imprimitivity \via\ homology}\label{s.homology}
As above, let $D\subset X$ be a reduced curve in a smooth projective
surface~$X$.
Denoting by $\inc\:D\into X$ the inclusion, let
\[*
\SS{D}=\Im[\inc_*\:H_2(D)\to H_2(X)/\Tors].
\]
As explained in \autoref{s.divisors},
$\SS{D}\subset\NS(X)$ is the subgroup generated by the
irreducible components of~$D$.
For convenience, we retain the notation $\inc\:D\into X$ and $\SS{D}$ in the
case when $D=\sum n_iD_i$, $n_i\ne0$, is a divisor in~$X$ (thus
identifying~$D$ with its support $\bigcup D_i$). The \emph{fundamental class} of
a divisor~$D$ is $[D]:=\sum n_i[D_i]$.


\theorem\label{th.reduction}
Let $\inc\:D\into X$ be as above, and assume that $H_1(X)=0$.
Then there are canonical isomorphisms
\[*
\Tors H_1(X\sminus D)=\Hom(\TT{D},\Q/\Z),\quad
H_1(X\sminus D)/\Tors=\Hom(\KK{D},\Z),
\]
where $\TT{D}:=\Tors(\NS(X)/\SS{D})$ and
$\KK{D}:=\Ker[\inc_*\:H_2(D)\to H_2(X)]$.
\endtheorem

\proof
By Poincar\'{e}--Lefschetz duality,
we have
$H_1(X\sminus D)=H^3(X,D)$.
Consider the following fragment of the cohomology
exact sequence of pair $(X,D)$:
\[*
H^2(X)\overset{\inc^*}\longto
H^2(D)\overset\Gd\longto
H^3(X,D)\longto
H^3(X).
\]
Since $H^3(X)=H_1(X)=0$ (Poincar\'{e} duality), we have
a canonical isomorphism
\[
H_1(X\sminus D)=\Coker\inc^*.
\label{eq.H1}
\]

As explained above, both $H^2(X)$ and $H^2(D)$ are
free abelian groups and, for both spaces, we have
$H^2(\cdot)=\Hom(H_2(\cdot),\Z)$; hence,
$\inc^*=\Hom(\inc_*,\id_\Z)$.
The exact sequence
\[*
0\longto\KK{D}\overset\inj\longto H_2(D)\overset{\inc_*}\longto H_2(X)
\]
can be regarded as a free resolution of $\QQ:=H_2(X)/\SS{D}$.
Applying $\Hom(\cdot,\Z)$, we obtain a cochain complex
\[*
0\longto H^2(X)\overset{\inc^*}\longto H^2(D)
 \overset{\inj^*}\longto\Hom(\KK{D},\Z)\longto0
\]
computing
the derived functors: $H^0=\Hom(\QQ,\Z)$,
$H^1=\Ext(\QQ,\Z)$, $H^i=0$ for $i\ge2$.
By the definition of $H^1$ and $H^2$, this gives us a short exact sequence
\[*
0\longto\Ext(\QQ,\Z)\longto\Coker\inc^*\longto\Hom(\KK{D},\Z)\longto0.
\]
Here, the first group is finite and the last one is free.
Hence,
\[*
\Ext(\QQ,\Z)=\Tors\Coker\inc^*\quad\text{and}\quad
\Hom(\KK{D},\Z)=\Coker\inc^*/\Tors.
\]
In view of~\eqref{eq.H1}, these two isomorphisms prove the two statements
of the theorem.
For the first statement, one should also observe that
$\Ext(\QQ,\Z)=\Ext(\Tors\QQ,\Z)$ (a property of finitely generated
abelian groups), $\Tors\QQ=\TT{D}$ (using the fact that $\NS(X)$ is
primitive in $H_2(S)$), and $\Ext(\TT{D},\Z)=\Hom(\TT{D},\Q/\Z)$
(apply $\Hom(\TT{D},\cdot)$ to the exact sequence $0\to\Z\to\Q\to\Q/\Z\to0$.)
\endproof

\subsection{The covering \pdfstr{F\sb{m}->F\sb1}{$\Fermat_m\to\Fermat_1$}}\label{s.m->1}
We make extensive use of the ramified covering
$\pr_m\:\Fermat_m\to\Fermat:=\Fermat_1$ given by
\[*
\pr_m\:[z_0:z_1:z_2:z_3]\mapsto[z_0^m:z_1^m:z_2^m:z_3^m].
\]
Clearly, $\Fermat$ is
the plane $\{z_0+z_1+z_2+z_3=0\}$, and $\pr_m$ is ramified over the union of
four lines $\rl_i:=\Fermat\cap\{z_i=0\}$, $i=0,1,2,3$.
The Galois group of~$\pr_m$ is $(\CG{m})^3$.
Assuming that $m\ge3$, the $3m^2$ lines in~$\Fermat_m$ are the irreducible
components of the preimage of the three lines
$\ln_i:=\Fermat\cap\{z_0+z_i=0\}$, $i=1,2,3$.
Introduce the divisors
$\ln:=\ln_1+\ln_2+\ln_3$,
$\rl:=\rl_0+\rl_1+\rl_2+\rl_3$,
and $\dv:=\ln+\rl$ in $\Fermat$.

With a further generalization in mind, redenote $\Fermat\ofm:=\Fermat_m$ and
consider the pull-backs $\ln_*\ofm:=\pr_m\1(\ln_*)$,
$\rl_*\ofm:=\pr_m\1(\rl_*)$, and $\dv\ofm:=\pr_m\1(\dv)$,
where $*$ is an appropriate subscript,
possibly empty.
Each\mnote{edited}
$\rl_j\ofm$ is a plane section
of $\Fermat\ofm$,
irreducible and reduced: it is the Fermat curve cut off
$\Fermat\ofm$ by the plane $\{z_j=0\}$.
On the other hand, $L\ofm$ also contains a number of plane sections, \eg,
those cut off by $\{z_i=\xi z_j\}$, $i\ne j$, $\xi^m=-1$.
Thus,
for any subset $J\subset\{0,1,2,3\}$, one has
\[
\SS{\dv\ofm}=\SS{\ln\ofm+\rl_J\ofm}=\SS{\ln\ofm}=\sset_m,
\label{eq.section}
\]
where $\rl_J\ofm:=\sum_{j\in J}\rl_j\ofm$.

Since $\rl$ is a
generic configuration of four lines in the plane~$\Fermat$, the
fundamental group
$\GG:=\pi_1(\Fermat\sminus\rl)$ equals $\Z^3$,
see~\cite[lemma in the proof of Theorem~8]{Zariski:group}.\mnote{ref
extended}
Since $\GG$ is abelian, from the Hurewicz theorem we have
$\GG=H_1(\Fermat\sminus\rl)=\Hom(\KK{\rl},\Z)$, see \autoref{th.reduction}.
This group
has four canonical generators~$\g_j$, $j=0,1,2,3$, \viz. the restrictions
to~$\KK{\rl}$ of the four generators of the group
$H^2(\rl)=\bigoplus_j\Z\cdot[\rl_j]^*$.
We have $\g_0+\g_1+\g_2+\g_3=0$, and $\GG$ is freely generated by $\g_1$,
$\g_2$, $\g_3$, \cf., \eg, \eqref{eq.H1}.\mnote{ref added}

An interesting generalization of the original question was suggested
in~\cite{Shimada.Takahashi:primitivity}.
Given an epimorphism $\Ga\:\GG\onto G$ to a finite abelian group $G$,
denote by $\pr\:\Fermat\of\Ga\to\Fermat$ the minimal resolution of
singularities of the ramified covering of~$\Fermat$ defined by~$\Ga$.
Let $\ln_*\of\Ga$, $\rl_*\of\Ga$, and $\dv\of\Ga$
be the pull-backs in~$\Fermat\of\Ga$ of $\ln_*$, $\rl_*$, and~$\dv$,
respectively. To be consistent with the previous notation, we regard an
integer~$m$ as the quotient projection $m\:\GG\onto\GG/m\GG$.
The components of $\dv\of\Ga$
(including the exceptional divisors)
represent some `obvious' elements of $\NS(\Fermat\of\Ga)$.
Using~\eqref{eq.Shioda} and the finite degree map
$\Fermat\ofm\dashrightarrow\Fermat\of\Ga$
defined by the inclusion $\Ker\Ga\subset m\GG$, $m:=\ls|G|$,
one has
\[
\SS{\dv\of\Ga}\otimes\Q=(\Pic\Fermat\of\Ga)\otimes\Q\quad
\text{whenever $\gcd(\ls|G|,6)=1$}.
\label{eq.Shioda.Ga}
\]
Thus, it is natural to ask whether $\SS{\dv\of\Ga}=\Pic\Fermat\of\Ga$,
or,
not assuming that $\ls|G|$ is prime to~$6$,
whether
$\SS{\dv\of\Ga}\subset\Pic\Fermat\of\Ga$ is a primitive subgroup.

\problem[Shimada--Takahashi~\cite{Shimada.Takahashi:primitivity}]\label{problem}
When does one have $\TT{\dv\of\Ga}=0$?
\endproblem

According to~\cite{Shimada.Takahashi:primitivity}, the answer to this
question is in the affirmative if the image~$G$ of~$\Ga$ is a cyclic group of
order $\ls|G|\le50$.
Another example is worked out in \autoref{s.toy}, see \autoref{th.toy}: the
answer is also in the affirmative if $\Ga(\g_i)=0$
for at least
one of the
standard generators $\g_i$, $i=0,1,2,3$.

\section{The Alexander module}\label{S.Am}

\subsection{The fundamental group}\label{s.pi1}
The\mnote{new remark}
line arrangement $\ln+\rl\subset\Cp2$ is well known; sometimes it is
referred to as \emph{Ceva-7}. Its fundamental group has been computed in many
ways and in many places; however, since we will work with a particular
presentation of this group, we repeat the computation here.

We will use the affine coordinates $x:=-z_1/z_0$, $y:=-z_3/z_0$
in the plane~$\Fermat$.
In these coordinates, $\rl_0$ becomes the line at infinity, and the other
components of $\dv$ are the lines of the form $\{r_xx+r_yy=r\}$
with $r_x,r_y,r\in\{0,1\}$,
see \autoref{fig.divisor}.
\figure
\centerline{\cpic{lines}}
\caption{The divisor $\dv:=\ln+\rl\subset\Fermat$}\label{fig.divisor}
\endfigure
The fundamental group $\pi_1:=\pi_1(\Fermat\sminus\dv)$ is easily
computed by the Zariski--van Kampen method~\cite{Zariski:group,vanKampen}.
Since we use a modified (or rather intermediate) version of
this approach, we outline briefly its proof, using $\dv$ as a
model.
(In full detail, the computation using the projection from a singular point is
explained, \eg, in~\cite{degt:book}.)\mnote{new citation}
Consider the projection
$p\:\Fermat\dashrightarrow\Cp1$, $(x,y)\mapsto x$.
This projection has four special fibers $F_a$, \viz. those over the points
$a\in\Delta:=\{-1,0,1,\infty\}$.
(Three of these fibers are components
of~$\dv$, but this fact is irrelevant for the moment.)
Let $F_*:=\bigcup F_a$, $a\in\Delta$.
Then the restriction
$p\:\Fermat\sminus(\dv\cup F_*)\to\Cp1\sminus\Delta$ is a locally
trivial fibration and,
since $\pi_2(\Cp1\sminus\Delta)=0$ and the fiber is connected,
Serre's exact sequence
(\emph{aka} long exact sequence of a fibration)\mnote{alternative name}
gives us a short exact sequence of fundamental groups
\[*
\{1\}\longto\pi_1(F\sminus\dv)\longto\pi_1(\Fermat\sminus(\dv\cup F_*))
 \longto\pi_1(\Cp1\sminus\Delta)\longto\{1\},
\]
where $F$ is a typical fiber of~$p$, \eg, the one over $x=\frac12$.
Choosing $(\frac12,-\frac32)$ for the basepoint, we have
$\pi_1(F\sminus\dv)=\<\vg_1,\vg_2,\vg_3,\vg_4\>$, see \autoref{fig.divisor}.
The group $\pi_1(\Cp1\sminus\Delta)$ is free, and the exact sequence splits.
A splitting can be constructed geometrically, identifying
$\pi_1(\Cp1\sminus\Delta)$ with $\pi_1(S\sminus F_*)=\<\hg_1,\hg_2,\hg_3\>$,
where $S\subset\Fermat$ is the section $y=-\frac32$,
the generators $\hg_1,\hg_2$ are as shown in
\autoref{fig.divisor}, and $\hg_3$ is a similar loop about
the fiber $F_{-1}$, not shown
in the figure. Thus, one arrives at the presentation
\[*
\pi_1(\Fermat\sminus(\dv\cup F_*))=
\bigl\<\vg_1,\vg_2,\vg_3,\vg_4,\hg_1,\hg_2,\hg_3\bigm|
 \hg_i\1\vg_j\hg_i=\Gb_i(\vg_j)\bigr\>,
\]
where $i=1,2,3$, $j=1,2,3,4$, and $\Gb_i\in\Aut\<\vg_1,\vg_2,\vg_3,\vg_4\>$
is the so-called \emph{braid monodromy}, \ie, the automorphism of the
fundamental group obtained by dragging the fiber along~$\hg_i$ while keeping
the basepoint in~$S$. (The formal definition is in terms of a trivialization
of the induced fibration
$(p\circ\hg_i)^*p$ over the segment $[0,1]$,
where $p\circ\hg_i$ is regarded as a map
$[0,1]\to\Cp1\sminus\Delta$; for all details,
see~\cite{Zariski:group,vanKampen}.)

Now, in order to pass to $\pi_1(\Fermat\sminus\dv)$, one needs to patch in
the only special fiber $F_{-1}$ that is not a component of~$\dv$.
This is done using the
Seifert--van Kampen theorem~\cite{vanKampen}.
In fact, the principal application of the theorem in~\cite{vanKampen} is the
following simple observation, which we state
in a slightly generalized form.

\lemma\label{lem.vK}
Let
$X$ be a smooth quasi-projective surface and $D\subset X$ a
closed smooth irreducible
curve. Then the inclusion homomorphism $\pi_1(X\sminus D)\to\pi_1(X)$ is an
epimorphism\rom; its kernel is normally generated by the class
$[\partial\Gamma]$, where $\Gamma$ is an analytic disc
transversal to~$D$ at its center and disjoint from~$D$ otherwise.
\done
\endlemma

Since $D$ is assumed irreducible, the conjugacy class of $[\partial\Gamma]$
in the statement does not depend on the choice of~$\Gamma$ or path connecting
$\partial\Gamma$ to the basepoint. The proof of the lemma is literally the
same as in~\cite{vanKampen}, using a tubular neighborhood of~$D$.

Applying \autoref{lem.vK}
to the curve $F_{-1}\sminus\dv$ in $\Fermat\sminus\dv$, we
obtain an extra relation $\hg_3=1$. In other words,
we disregard the generator
$\hg_3$ and convert the four relations $\hg_3\1\vg_j\hg_3=\Gb_3(\vg_j)$ into
$\vg_j=\Gb_3(\vg_j)$, $j=1,2,3,4$.

The computation of the braid monodromy is straightforward and well known,
\eg, using equations of the lines; it is left to the reader. (Essentially,
it is
the braid monodromy of the nodal arrangement
$\ln_1+\ln_2+\rl_2+\rl_3$ of four lines.) Denoting
by $\Gs_1,\Gs_2,\Gs_3$
the Artin generators~\cite{Artin} of the braid group $\BG4$ acting on
$\<\vg_1,\vg_2,\vg_3,\vg_4\>$, we have
$\Gb_1=\Gs_1^2\Gs_3^2$,
$\Gb_2=\Gs_2^2$,
and $\Gb_3=\Gs_1\1\Gs_3\1\Gs_2^2\Gs_3\Gs_1$.\mnote{typos corrected}
(It is worth recalling that, assuming the left action of the automorphism
group, the braid monodromy is an \emph{anti}-homomorphism
$\pi_1(\Cp1\sminus\Delta)\to\BG4$.)
Indeed,\mnote{a few more words about the computation}
$\Gb_1$ and $\Gb_2$ are essentially computed in
the very first paper on the subject, \viz.~\cite{Zariski:group}:
each is the local monodromy about a simple node (one full twist of a pair of
points about their barycenter) or a pair of disjoint nodes.
The remaining braid~$\Gb_3$ is the local monodromy about another node, which
is
translated to the common reference fiber along the real axis;
when circumventing the singular fiber at the origin, it gets conjugated
by `one half' of~$\Gb_1$, which is $\Gs_1\Gs_3$.

Putting everything together,\mnote{par break; introductory words changed}
after a slight simplification the
nontrivial relations for the fundamental group $\pi_1(\Fermat\sminus\dv)$
take the form
\begin{gather}
[\hg_2,\vg_1]=[\hg_2,\vg_4]=1,\label{rel.1}\\
\hg_2\vg_2\vg_3=\vg_2\vg_3\hg_2=\vg_3\hg_2\vg_2\label{rel.2}
\end{gather}
(the relations $\hg_2\1\vg_j\hg_2=\Gb_2(\vg_j)$ from the fiber $x=1$),
\begin{gather}
\hg_1\vg_1\vg_2=\vg_1\vg_2\hg_1=\vg_2\hg_1\vg_1\rlap,\label{rel.3}\\
\hg_1\vg_3\vg_4=\vg_3\vg_4\hg_1=\vg_4\hg_1\vg_3\label{rel.4}
\end{gather}
(the relations $\hg_1\1\vg_j\hg_1=\Gb_1(\vg_j)$ from the fiber $x=0$), and
\[
[\vg_2^{-1}\vg_1\vg_2,\vg_4]=1\label{rel.5}
\]
(the relations $\vg_j=\Gb_3(\vg_j)$ from the fiber $x=-1$).

By \autoref{lem.vK},
the inclusion $\inj\:\Fermat\sminus\dv\into\Fermat\sminus\rl$ induces
the map
\[*
\inj_*\:\pi_1\onto\GG:\qquad
\hg_1\mapsto\g_1,\quad
\vg_2\mapsto\g_2,\quad
\vg_3\mapsto\g_3,\quad
\hg_2,\vg_1,\vg_4\mapsto0.
\]

\subsection{The `universal' covering}\label{s.covering}
Throughout\mnote{par added}
the paper we use freely the following well-known fact, often
referred to as \emph{theory of covering spaces}: for any connected,
locally path connected, and micro-simply connected topological space~$X$
(\eg, for any connected simplicial complex) with a basepoint $x_0\in X$,
there is a natural
equivalence between the category of coverings
$(\tilde X,\tilde x_0)\to(X,x_0)$ and covering maps
(identical on $X$) and
that of subgroups of $\pi_1(X,x_0)$ and inclusions.
If the subgroup is normal (\emph{regular}, or \emph{Galois} coverings), it
can be described as the kernel of an epimorphism $\Ga\:\pi_1(X,x_0)\onto G$;
the image~$G$ serves then as the group of the deck translations of the
covering.

Consider
an epimorphism $\Ga\:\GG\onto G$. In this section, we do not assume
$G$ finite; in fact, we start with a study of the `universal'
$\GG$-covering, corresponding to the identity map $0\:\GG\onto\GG/0\GG=\GG$.
(Admittedly awkward,\mnote{an apology added}
this notation is compliant with $m\:\GG\onto\GG/m\GG$
introduced earlier.)


Consider the composition
\[*
\bGa\:\pi_1\overset{\inj_*}{\relbar\joinrel\onto}
\pi_1(\Fermat\sminus\rl)=\GG\overset\Ga{\relbar\joinrel\onto}G
\]
and denote by
$\Fermat^\circ\of\Ga$ the $G$-covering of $\Fermat\sminus\dv$ defined
by~$\bGa$. By the Hurewicz theorem, $H_1(\Fermat^\circ\of\Ga)$ is
the abelianization of
$\pi_1(\Fermat^\circ\of\Ga)=\Ker\bGa$.
The action of the deck translations of the covering makes this group
a $\Z[G]$-module;
regarded as such,
it is often referred to as the \emph{Alexander module} of~$\bGa$.

The construction of the Alexander module
fits into a more general framework and admits a purely algebraic description.
Consider a group $\pi$ and an epimorphism $\bGa\:\pi\onto G$ to an
abelian group~$G$. Then the Alexander module of~$\bGa$ is the abelian
group $\Am:=\Ker\bGa/[\Ker\bGa,\Ker\bGa]$ regarded as a $\Z[G]$-module \via\
the
$G$-action
defined as follows: given $a\in\Am$ and $g\in G$, the image $g(a)$ is
the class in $\Am$ of the element
$\tilde g\tilde a\tilde g\1\in\Ker\bGa$, where
$\tilde a,\tilde g\in\pi$ are some lifts of $a,g$, respectively.
This class does not depend on the choice of the lifts, and the action is well
defined.

Crucial is the fact that $H_1(\Fermat^\circ\of\Ga)$ depends on the epimorphism
$\bGa\:\pi_1\onto G$ only. Hence, we can replace
$\Fermat\sminus\dv$
with any \CW-complex~$X$ with $\pi_1(X)=\pi_1$.
We take for~$X$ a space with a single $0$-cell~$e^0$, one
$1$-cell $e^1_i\in\{\ag_1,\ag_2,\ag_3,\cg_1,\cg_2,\cg_3\}$
for each of the six generators $\hg_1,\vg_2,\vg_3,\hg_2,\vg_4,\vg_1$
of~$\pi_1$ (in the order listed),
and one $2$-cell~$e^2_j$ for each
relation~\eqref{rel.1}--\eqref{rel.5}.
In the $\GG$-covering~$X\of0$, each cell~$e$ gives rise to a whole
$\GG$-orbit $\{g\otimes e\,|\,g\in\GG\}$.
(For the moment,
the symbols $g\otimes e$ are merely cell labels; we only
assume that the labelling is compatible with the $\GG$-action, \ie, for any
cell~$e$ in~$X$ and pair $h,g\in\GG$ we have
$h(g\otimes e)=(h+g)\otimes e$.)

Following the tradition,
let us
identify $\Z[\GG]$ with the ring
\[*
\GL:=\Z[t_1^{\pm1},t_2^{\pm1},t_3^{\pm1}]
\]
of Laurent polynomials,
where
the variables $t_1$, $t_2$, $t_3$ correspond to the generators
$\hg_1\mapsto\g_1$, $\vg_2\mapsto\g_2$, $\vg_3\mapsto\g_3$
about
$\rl_1$, $\rl_2$, $\rl_3$,
respectively.
In other words, we identify~$\GG$ with the multiplicative abelian group
generated by $t_1,t_2,t_3$; we will also use this multiplicative notation in
the cell labels.
We can assume, in addition, that the
labelling is chosen so that
the left end of each `initial' $1$-cell $1\otimes e$
is attached to $1\otimes e^0$, \ie,
$(1\otimes e)(0)=1\otimes e^0$. (Here, we regard an oriented $1$-cell as a
path $[0,1]\to X\of0$.)
Then,
from the definition of the covering it follows that
the right ends are attached as follows:
\[
(1\otimes\ag_i)(1)=t_i\otimes e^0,\quad
(1\otimes\cg_j)(1)=1\otimes e^0,\quad i,j=1,2,3,
\label{eq.ends}
\]
\ie, the generators $\hg_1,\vg_2,\vg_3$ are `unwrapped', whereas
$\hg_2,\vg_1,\vg_4$ remain `latent'.
The other ends are determined by the $\GG$-action: for a $1$-cell~$e$ in~$X$,
a monomial $t$ in $t_1,t_2,t_3$, and $\epsilon=0,1$ we have
$(t\otimes e)(\epsilon)=t((1\otimes e)(\epsilon))$.

Recall that the member~$C_n$ of the cellular chain complex associated to a
\CW-complex $Y$ is the free abelian group generated by the $n$-cells of~$Y$.
Thus, each cell $e$ of~$X$ gives rise to a direct summand
$\bigoplus\Z(g\otimes e)$, $g\in\GG$, in the complex of $X\of0$; this summand is
naturally identified with the free $\GL$-module $\GL e$.
(It is this identification that explains the usage of $\otimes$ in the
labels.)
Furthermore, since the \CW-structure on $X\of0$ is $\GG$-invariant, the
boundary homomorphisms are $\GL$-linear.
Thus, the
chain complex $C_*:=C_*\of0$ of
$X\of0$ is a
complex of free $\GL$-modules
of the form
\[*
0\longto
C_2\overset{\partial_2}\longto
C_1=\GL\ag_1\oplus\GL\ag_2\oplus\GL\ag_3\oplus\GL\cg_1\oplus\GL\cg_2\oplus\GL\cg_3
 \overset{\partial_1}\longto
C_0=\GL\longto0
\]
(we omit the generator $e^0$ of $C_0$),
where $\partial_1$ is given by~\eqref{eq.ends}:
\[
\partial_1\ag_i=(t_i-1),\quad
\partial_1\cg_j=0,\quad i,j=1,2,3.
\label{eq.d1}
\]
The module~$C_2$ has nine generators, of which six have non-trivial
images under~$\partial_2$:
\begin{gather}
(t_2t_3-1)\cg_1,\label{rel.b1}\\
 (t_3-1)\cg_1+(t_3-1)\ag_2-(t_2-1)\ag_3\label{rel.a1}
\end{gather}
from~\eqref{rel.2},
\begin{gather}
(t_1t_3-1)\cg_2,\label{rel.b2}\\
 (t_3-1)\cg_2+(t_3-1)\ag_1-(t_1-1)\ag_3\label{rel.a2}
\end{gather}
from~\eqref{rel.4}, and
\begin{gather}
(t_1t_2-1)\cg_3,\label{rel.b3}\\
 (t_1-1)\cg_3+(t_1-1)\ag_2-(t_2-1)\ag_1\label{rel.a3}
\end{gather}
from~\eqref{rel.3}.
Relations~\eqref{rel.1} and~\eqref{rel.5} contribute~$0$ to~$\Im\partial_2$.

\example\label{ex.relations}
The proof
of \eqref{rel.b1}--\eqref{rel.a3} is a straightforward computation.
As an example, consider~\eqref{rel.2}, which can be written in
the form of two relations
\[*
\hg_2\vg_2\vg_3\hg_2\1\vg_3\1\vg_2\1=1,\quad
\hg_2\vg_2\vg_3\vg_2\1\hg_2\1\vg_3\1=1.
\]
The word in the left hand side
of the first relation
corresponds to the sequence
$\cg_1$, $\ag_2$, $\ag_3$, $\cg_1\1$, $\ag_3\1$, $\ag_2\1$
of $1$-cells in~$X$ along which a
$2$-cell~$e^2_1$ is attached.
(The inverse for a $1$-cell means the reversion of the orientation.)
Lift this sequence to $X\of0$,
\emph{starting each cell at the end of the previous one},
see~\eqref{eq.ends}:
\[*
1\otimes\cg_1,\
1\otimes\ag_2,\
t_2\otimes\ag_3,\
(t_2t_3\otimes\cg_1)\1,\
(t_2\otimes\ag_3)\1,\
(1\otimes\ag_2)\1.
\]
(Observe that,
for example, $t_2\otimes\ag_3$ connects $t_2\otimes e^0$ to $t_2t_3\otimes
e^0$, see~\eqref{eq.ends};
hence, the lift of $\ag_3\1$ \emph{starting} at $t_2t_3\otimes e^0$ is
$(t_2\otimes\ag_3)\1$; it \emph{ends} at $t_2\otimes e^0$.
Note also that $(1\otimes\ag_2)\1$ ends at $1\otimes e^0$, \ie, the lift is a
loop, as expected.)
We obtain a sequence of $1$-cells along which a $2$-cell in $X\of0$, \viz.
one of the lifts of~$e^2_1$,
is
attached; writing this sequence as a chain, we get
$\partial_2e^2_1=(1-t_2t_3)\cg_1\in C_1$, which is~\eqref{rel.b1} up to sign.
Similarly, the second relation lifts to
the sequence
\[*
1\otimes\cg_1,\
1\otimes\ag_2,\
t_2\otimes\ag_3,\
(t_3\otimes\ag_2)\1,\
(t_3\otimes\cg_1)\1,\
(1\otimes\ag_3)\1,
\]
which gives us~\eqref{rel.a1}.
\endexample

\subsection{Other coverings}\label{s.others}
Now,
given an epimorphism $\Ga\:\GG\onto G$,
it induces a ring homomorphism
$\Ga_*\:\GL\onto\Z[G]$, making $\Z[G]$ a
$\GL$-module. Clearly, the $G$-covering $X\of\Ga$ is the quotient
space $X\of0/\!\Ker\Ga$, the cells in $X\of\Ga$ being the $\Ker\Ga$-orbits of
those in $X\of0$.
The chain homomorphism $C_*\to C_*(X\of\Ga)$ induced by the quotient
projection merely identifies
the basis elements
(which are the cells)
within each orbit of $\Ker\Ga$;
algebraically, it can be expressed
as
the tensor product
\[*
{\id}\otimes\Ga_*\:
C_*=C_*\otimes_\GL\GL
 \longto
 C_*\otimes_\GL\Z[G]=C_*(X\of\Ga).
\]

Recall, see the beginning of \autoref{s.covering}, that the $1$-homology of
the covering spaces depend only on the homomorphism
$\bGa\:\pi_1\onto G$.
Hence, the group $H_1(\Fermat^\circ\of\Ga)=H_1(X\of\Ga)$ is computed by the
complex
$C_*\of\Ga:=C_*\otimes_\GL\Z[G]$.
In view of the right exactness
$\Coker(\partial_2\otimes_\GL\Ga_*)=(\Coker\partial_2)\otimes_\GL\Z[G]$,
our primary interest is the quotient
$\Hm\of\Ga:=C_1\of\Ga/\Im\partial_2$. Explicitly, $\Hm\of\Ga$
can be described as the $\GL$-module generated by the six elements
$\ag_1,\ag_2,\ag_3,\cg_1,\cg_2,\cg_3$
that are
subject to relations~\eqref{rel.b1}--\eqref{rel.a3}
and the extra relation
\[
t_1^{r_1}t_2^{r_2}t_3^{r_3}=1\quad
\text{whenever $\Ga(r_1\g_1+r_2\g_2+r_3\g_3)=0$}.
\label{rel.Ga}
\]
Summarizing,
after the identification $C_0\of\Ga=\Z[G]$ and $H_0(X\of\Ga)=\Z$, we have an
exact sequence
\[
0\longto
 H_1(\Fermat^\circ\of\Ga)\longinto
 \Hm\of\Ga\overset{\partial_1}\longto
 \Z[G]\longto
 \Z\longto0,
\label{eq.H1.final}
\]
where the last homomorphism is the augmentation $g\mapsto1$, $g\in G$.

Recall that the \emph{rank} $\rank A$ of a finitely generated abelian
group~$A$ is the maximal number of linearly independent elements of~$A$,
whereas its \emph{length} $\ell(A)$ is the minimal
number of elements generating~$A$. One has $\rank A=\ell(A)$ if and only if
$A$ is free.


\lemma\label{lem.H}
For any epimorphism $\Ga\:\GG\onto G$, there is a natural isomorphism
$\Tors H_1(\Fermat^\circ\of\Ga)=\Tors\Hm\of\Ga$.
If $G$ is finite, then
$\ell(H_1(\Fermat^\circ\of\Ga))=\ell(\Hm\of\Ga)-\ls|G|+1$
and $\rank H_1(\Fermat^\circ\of\Ga)=\rank\Hm\of\Ga-\ls|G|+1$.
\endlemma

\proof
Since $\Im\partial_1\subset\Z[G]$ is a free abelian
group, the inclusion in~\eqref{eq.H1.final} induces an isomorphism of the
torsion parts.
This isomorphism and the obvious fact that
$\ell(A)=\rank A+\ell(\Tors A)$ for any
finitely generated abelian group~$A$
imply that the length
and rank identities in the statement are equivalent to each other.
The rank identity follows from the additivity of rank in~\eqref{eq.H1.final}
and the observation that $\rank\Z[G]=\ls|G|$.
\endproof

%

\subsection{Fermat surfaces}\label{s.Fermat}
If the image~$G$ of $\Ga\:\GG\onto G$ is finite, one obviously has
$\Fermat^\circ\of\Ga=\Fermat\of\Ga\sminus\dv\of\Ga$.
If $\Ga=m\in\Z$, \ie, in the case of a classical Fermat surface $\Fermat\ofm$,
it is more convenient to consider a smaller divisor
$\tln\ofm:=\ln\ofm+\rl_0\ofm$, see~\eqref{eq.section}.\mnote{ref added}
The
fundamental
group $\pi_1(\Fermat\ofm\sminus\tln\ofm)$ is given by
\autoref{lem.vK}:
it is the quotient of $\Ker\bGa=\pi_1(\Fermat^\circ\of\Ga)$
by the extra relations
$\hg_1^m=\vg_2^m=\vg_3^m=1$ (as\mnote{explanation added}
the ramification index at each
component of $\rl\ofm$ is obviously~$m$).
Hence, the homology group
$H_1(\Fermat\ofm\sminus\tln\ofm)$ can be computed using the complex $C_*\ofm$
with three extra $2$-cells $e^2_i$, $i=1,2,3$, mapped by $\partial_2$ to
$\poly_m(t_i)\ag_i$, where
\[*
\poly_n(t):=(t^n-1)/(t-1),\quad n\in\Z.
\]
This computation is similar to \autoref{ex.relations}:
for example, the loop $\hg_1^m$ lifts to the sequence
$1\otimes\ag_1, t_1\otimes\ag_1,t_1^2\otimes\ag_1,\ldots,t_1^{m-1}\otimes\ag_1$
of $1$-cells,
which results in the chain
$(1+t_1+t_1^2+\ldots+t_1^{m-1})\ag_1=\poly_m(t_1)\ag_1\in C_1\ofm$.
Note that this chain is a cycle, as in $C_1\ofm$ we have the relation
$t_1^m=1$.

\remark\label{rem.Z->GL}
Strictly speaking,
the new complex is that of abelian groups rather than $\GL$-modules, as we
add three $2$-cells only, \ie, three summands $\Z e^2_i$ in $C_2\ofm$.
However, in the presence of the relations $t_i^m=1$, $i=1,2,3$,
\cf.~\eqref{rel.Ga},
one can use \eqref{rel.b1}--\eqref{rel.a3} to show that all three images
$\poly_m(t_i)\ag_i$ are $\GG$-invariant.
Hence, without changing the $1$-homology of the complex, we can formally
replace each summand $\Z e^2_i$ with $\GL e^2_i$, extending~$\partial_2$ by
$\GL$-linearity.
Geometrically, we replace a single disk $\Gamma$ as in
\autoref{lem.vK} with a $G$-orbit consisting of $m^3$ disks.
Since the curve~$\rl_i\ofm$ patched in is irreducible
(all disks intersecting the same component),
this change does not affect
the fundamental group.
\endremark

Now, as in \autoref{s.others}, instead of extending the $C_2$-term
of the complex, we can
add extra relations to~$C_1$.
Summarizing, we have
\[*
H_1(\Fermat\ofm\sminus\tln\ofm)=\Ker[\partial_1\:\tAm\ofm\to C_0\ofm],
\]
where $\tAm\ofm$ is the quotient of $\Am\of0$ by the extra relations
\[
t_i^m=1,\quad
\poly_m(t_i)\ag_i=0,\quad
i=1,2,3.
\label{rel.m}
\]
Arguing as in the proof of \autoref{lem.H}, we obtain the
identity
\[
\ell(H_1(\Fermat\ofm\sminus\tln\ofm))=\ell(\tAm\ofm)-m^3+1.
\label{eq.H}
\]

\subsection{Other Delsarte surfaces}\label{s.finite}
In the generalized case,
the first question that arises is whether \autoref{th.reduction} is
applicable, \ie, whether $H_1(\Fermat\of\Ga)=0$.
To state the result, introduce the following notation:
given a pair of integers $0\le i,j\le3$, let $\GG_{ij}:=\Z\g_i\oplus\Z\g_j\subset\GG$,
where $\g_i\in\GG$ are the canonical generators, see \autoref{s.m->1}.

Recall
that the blow-up $\Gs\:\tilde X\to X$ of a \emph{smooth} point of a
surface~$X$ induces an isomorphism of both the fundamental group~$\pi_1$ and
first homology group~$H_1$ of the surface. Hence, up to canonical
isomorphism, the groups $\pi_1$ and $H_1$ do not depend on the resolution of
singularities.

\proposition\label{prop.H1}
For an epimorphism $\Ga\:\GG\onto G$, $\ls|G|<\infty$, one
has
\[*
\textstyle
\pi_1(\Fermat\of\Ga)=
H_1(\Fermat\of\Ga)=\Ker\Ga/\sum\GG_{ij}\cap\Ker\Ga,
\]
the summation running over all pairs $0\le i,j\le3$ of integers.
\endproposition

\proof
We start with the abelian group
$\pi_1(\Fermat\sminus\rl)=\GG$ generated by
$\hg_1\mapsto\g_1$, $\vg_2\mapsto\g_2$, $\vg_3\mapsto\g_3$,
see \autoref{s.pi1}.
Clearly,
$\pi_1(\Fermat\of\Ga\sminus\rl\of\Ga)=H_1(\Fermat\of\Ga\sminus\rl\of\Ga)=\Ker\Ga$.
(This group can also be regarded as a $\GL$-module, but the module structure
is trivial: $t_1=t_2=t_3=1$.) For the rest of the proof, we use the
additive notation for the fundamental group (as the groups of interest
are subquotients
of~$\GG$).

Let $\Fermat'\of\Ga$ be the manifold obtained from
$\Fermat\of\Ga\sminus\rl\of\Ga$ by patching the components of the
proper transform of $\rl\of\Ga$ \emph{away from the exceptional divisor}.
At a generic point of~$\rl_i$,
the ramification index~$m_i$
of the ramified covering $\Fermat\of\Ga\to\Fermat$
equals the index $[\GG_{ii}:\GG_{ii}\cap\Ker\Ga]$,
$i=0,1,2,3$.
Hence, by \autoref{lem.vK}, the inclusion induces
an epimorphism $\Ker\Ga\onto\pi_1(\Fermat'\of\Ga)$ whose kernel is
generated by the elements $m_i\g_i$. Thus, we have an
isomorphism
\[
\textstyle
\pi_1(\Fermat'\of\Ga)=\Ker\Ga/\sum_i\GG_{ii}\cap\Ker\Ga,\quad
i=0,1,2,3.
\label{eq.pi1F}
\]
(Strictly speaking,
unlike the case of the Fermat surfaces,
the curve $\rl_i\of\Ga$ may
be reducible, so that we need to
attach a separate disk~$\Gamma$ as in \autoref{lem.vK} for each component
of this curve.
However, since the $G$-action is trivial in the $1$-homology
$H_1=\pi_1$, all disks
result in the same relation $m_i\g_i=0$, \cf. \autoref{rem.Z->GL}.)

What remains is patching the exceptional divisors.
Fix a pair $0\le i<j\le3$ and
let $\tS$ be a singular point of the normalized, but yet unresolved ramified
covering over the point $S:=\rl_i\cap\rl_j$.
Fix a resolution of singularities and
let~$E$ be the exceptional divisor over~$\tS$.
Pick a sufficiently small ball $U\subset\Fermat$ about~$S$
and denote by~$\tU$ the connected component of the preimage
of~$U$ containing~$E$.
With respect to an appropriate smooth triangulation,
$\tU$ is a regular neighborhood
of~$E$; hence, $E$ is a strict deformation retract of~$\tU$, $\tU\sim E$.
On the other
hand, $\tU$ is a $4$-manifold with boundary $\partial\tU$,
and the latter is a
covering of the $3$-sphere $\partial U$ ramified over the Hopf link
$\rl\cap\partial U$.

Note also that the contraction of~$E$ gives us the space $\tU/E$
which is the cone over~$\partial\tU$ (with the
vertex $\tS=E/E$); hence, we have a homotopy equivalence
(strict deformation retraction)
$\tU\sminus E=(\tU/E)\sminus\tS\sim\partial\tU$.

We have $\pi_1(\partial U\sminus\rl)=\GG_{ij}$ and,
hence, $\pi_1(\partial\tU\sminus\rl\of\Ga)=\GG_{ij}\cap\Ker\Ga$.
As above, similar to \autoref{lem.vK}, patching the union of circles
$\partial\tU\cap\rl\of\Ga$ results in the pair of relations
$m_i\g_i=m_j\g_j=0$. Thus,
\[
\pi_1(\partial\tU)=(\GG_{ij}\cap\Ker\Ga)/(\GG_{ii}\cap\Ker\Ga+\GG_{jj}\cap\Ker\Ga)
\label{eq.pi1dU}
\]
is a finite group. Then $H_1(\partial\tU;\Q)=0$, \ie,
$\partial\tU$ is a rational
homology sphere and
$\tS$ is a rational singular point.
For us, important is the fact that
$\pi_1(\tU)=\pi_1(E)=0$, which can easily be proved
directly. Indeed, since $\tU\sim E$ and $\dim_\R E=2$, we have
$H^3(\tU;\Q)=0$; then also $H_1(\tU,\partial\tU;\Q)=0$ (Lefschetz duality),
and the fragment
\[*
H_1(\partial\tU;\Q)\longto H_1(\tU;\Q)\longto H_1(\tU,\partial\tU;\Q)
\]
of the homology exact sequence of pair $(\tU,\partial\tU)$ implies
$H_1(\tU;\Q)=H_1(E;\Q)=0$.
On the other hand, $E$ is a
connected projective algebraic curve, and it is easily seen
that $E$ is homotopy equivalent to the wedge of closed topological surfaces
(the components of the normalization of~$E$) and a number of
circles.
(Roughly, we can `blow-up' the locally reducible singular points of~$E$ to
line segments, separating the analytic branches and replacing~$E$ with a
disjoint union of topologically nonsingular closed surfaces with a number of
segments attached.
Then, within each surface, move the ends of the segments to a single point.
Finally, contract several segments to make the surfaces share a
common point;
the result is a wedge as stated.)
For such a wedge $E\sim\bigvee_iE_i$, all groups are easily computed (\eg,
using
iteratedly  the Mayer--Vietoris exact sequence~$(8.8)$ in \cite[Chapter III]{Dold}
and Seifert--van Kampen theorem~\cite{vanKampen}, or just decomposing
the wedge into cells):
\[*
\textstyle
H_1(E;\Q)=\bigoplus_iH_1(E_i;\Q),\quad
\pi_1(E)=\mathop*_i\pi_1(E_i).
\]
Clearly, $H_1(E;\Q)=0$
if and only if all surface components are $2$-spheres and there are no circles
present. Then obviously $\pi_1(E)=0$.

Now, start with $\Fermat'\of\Ga$
and proceed patching the
exceptional divisors one by one. Let $\Fermat''$ be an intermediate space,
not yet containing~$E$.
Applying the Seifert--van Kampen theorem~\cite{vanKampen} to the union
$\Fermat''\cup\tU$ and using the homotopy equivalence
$\Fermat''\cap\tU=\tU\sminus E\sim\partial\tU$,
we obtain the amalgamated free product
\[*
\pi_1(\Fermat''\cup\tU)
 =\bigl(\pi_1(\Fermat'')*\pi_1(\tU)\bigr)/\pi_1(\partial\tU)
 =\pi_1(\Fermat'')/(\GG_{ij}\cap\Ker\Ga).
\]
(For the last isomorphism, we use~\eqref{eq.pi1dU} and
the identity
$\pi_1(\tU)=0$.)
The group $\pi_1(\Fermat'\of\Ga)$ is given by~\eqref{eq.pi1F} and,
after all the exceptional divisors have been patched, we arrive at the
expression in the statement.
\endproof

If $H_1(\Fermat\of\Ga)=0$, \autoref{th.reduction} and \autoref{lem.H} imply
that
\[
\TT{\dv\of\Ga}\cong\Tors\Hm\of\Ga.
\label{eq.T=Tors}
\]
Unfortunately, as a $\Z[G]$-module, $\Hm\of\Ga$ is far from free
and it is difficult to control its $\Z$-torsion.
(Experiments show that, at least, the intermediate quotients similar to those
considered in \autoref{lem.ell} do often have torsion.)
An attempt of a direct computation is made in~\autoref{s.toy}, whereas in the
case of the classical Fermat surfaces we have to take a detour and estimate the
length instead.
The following two exact sequences may prove useful:
\[*
\Hm\of\Ga\overset{\partial_1}\longto
\Z[G]\overset{\epsilon}\longto\Z\longto0,
\]
where $\epsilon$ is the augmentation, see~\eqref{eq.H1.final},
and\mnote{edited}
\[*
0\longto\Bm\of\Ga\longto\Ker\partial_1\longto\Ker\Ga\longto0,
\]
where $\Bm\of\Ga\subset\Hm\of\Ga$ is the submodule generated by~$\cg_1$,
$\cg_2$,~$\cg_3$.
The\mnote{explanation added}
former sequence merely states that
$H_0(C_*\of\Ga)=H_0(\Fermat^\circ\of\Ga)=\Z$.
For the latter, we patch $\ln\of\Ga$ (by using \autoref{lem.vK} or merely
forgetting the generators $\hg_2,\vg_1,\vg_4$, hence $\cg_1,\cg_2,\cg_3$ in
the first place) to compute the group
$H_1(\Fermat\of\Ga\sminus\rl\of\Ga)=\pi_1(\Fermat\of\Ga\sminus\rl\of\Ga)=\Ker\Ga$;
the resulting complex is $0\to\Hm\of\Ga/\Bm\of\Ga\to\Z[G]\to0$.
Both sequences split, and we can extend~\eqref{eq.T=Tors}
to
\[
\TT{\dv\of\Ga}\cong\Tors\Bm\of\Ga=\Tors\Hm\of\Ga,
\label{eq.T=Tors0}
\]
still under the assumption that $H_1(\Fermat\of\Ga)=0$.

\section{Proof of the main theorem}\label{S.proof}

\subsection{The length of \pdfstr{A[m]}{$\tAm\ofm$}}\label{s.length}
Fix an integer $m\ge1$
and
consider the $\GL$-module $\tAm\ofm$
introduced in \autoref{s.Fermat}.
Recall that $\tAm\ofm$ is generated by
six elements $\ag_i$,
$\cg_j$, $i,j=1,2,3$, subject to the relations \eqref{rel.b1}--\eqref{rel.a3}
and~\eqref{rel.m}.
Observe that relations \eqref{rel.b1}, \eqref{rel.b2}, and~\eqref{rel.b3} can
be recast in the form
\[
t_i\cg_k=t_j\1\cg_k\quad
\text{whenever $\{i,j,k\}=\{1,2,3\}$}.
\label{rel.bb}
\]

We introduce a few \latin{ad hoc} notations. Given $i=1,2,3$, let
\[*
\GL_i:=\Z[t_i]/(t_i^m-1),\qquad
\tGL_i:=\Z[t_i]/\poly_m(t_i).
\]
These rings can be regarded as $\GL$-modules, but we usually do not specify
the action of the other two variables: it varies from case to case.
In fact, we repeatedly use the following simple observation, which is an
immediate consequence of~\eqref{rel.bb}.

\lemma\label{lem.b-module}
Let $i,j,k\in\{1,2,3\}$, $k\ne i$, and $p\in\GL$, and let $\Am$ be a
subquotient of $\tAm\ofm$ generated by a single element $x:=p\cg_i$. Assume
that
either $t_j=1$ or $t_i=t_k^{\pm1}$ on~$\Am$.
Then $\Am$ is a quotient of $\GL_sx$ for an appropriate index
$s\in\{1,2,3\}$.

If $x$ is also annihilated by
$\poly_m(t_s)$, then $\Am$ is a quotient of $\tGL_sx$.
\done
\endlemma

The precise description of the `appropriate' index~$s$
(not necessarily unique) is left to the reader.
Clearly, $\ell(\GL_s)=m$ and $\ell(\tGL_s)=m-1$.

For a generator
$x\in\{\ag_1,\ag_2,\ag_3,\cg_1,\cg_2,\cg_3\}$,
let
\[*
x':=(t_1-1)x,\quad
\3x:=(t_3-1)x,\quad
\3x':=(t_1-1)\3x.
\]
Observe that always
\[
\poly_m(t_1)x'=\poly_m(t_3)\3x=\poly_m(t_1)\3x'=\poly_m(t_3)\3x'=0.
\label{rel.x}
\]

We will use a filtration
$0=\Am_0\subset\Am_1\subset\ldots\subset\Am_7=\tAm\ofm$,
where $\Am_k\subset\tAm\ofm$ are the submodules
defined in \autoref{lem.ell} below.

Let $\Gd_m:=1$ if $m$ is even and $\Gd_m:=0$ if $m$ is odd.

\lemma\label{lem.ell}
One has the following equations and inequalities\rom:
\roster
\item\label{ell.1}
$\ell(\Am_1/\Am_0)=m^3-m^2$,
where $\Am_1$ is the submodule generated by~$\ag_3$\rom;
\item\label{ell.2}
$\ell(\Am_2/\Am_1)\le3(m-1)-\Gd_m$,
where $\Am_2:=\Am_1+\GL\bag_2'+\GL\bbg_3'$\rom;
\item\label{ell.3}
$\ell(\Am_3/\Am_2)\le3(m-1)$,
where $\Am_3:=\Am_1+(t_3-1)\tAm\ofm$\rom;
\item\label{ell.4}
$\ell(\Am_4/\Am_3)=m^2-m$,
where $\Am_4:=\Am_3+\GL\ag_1$\rom;
\item\label{ell.5}
$\ell(\Am_5/\Am_4)\le m-1$,
where $\Am_5:=\Am_4+\GL\ag_2'+\GL\cg_3'$\rom;
\item\label{ell.6}
$\ell(\Am_6/\Am_5)=m-1$,
where $\Am_6:=\Am_5+\GL\ag_2$\rom;
\item\label{ell.7}
$\ell(\Am_7/\Am_6)\le2m+1$,
where $\Am_7:=\tAm\ofm$.
\endroster
Hence, $\ell(\Am)\le m^3+9m-7-\Gd_m$.
\endlemma

\proof
One has $\ell(\Am_1)\le m^2(m-1)$ due to~\eqref{rel.m}. On the other hand, the
boundary homomorphism $\partial_1$
maps $\Am_1$ onto
$(t_3-1)C_0\ofm$. Hence, there are no other relations
in~$\Am_1$, and Statement~\iref{ell.1} holds.
Furthermore, $\partial_1$ factors to a homomorphism
\[*
\tAm\ofm/\Am_3\to C_0':=C_0\ofm/(t_3-1)
\]
which maps $\Am_4/\Am_3$
isomorphically onto $(t_1-1)C_0'$, proving Statement~\iref{ell.4}.
Then, $\partial_1$ factors to
\[*
\tAm\ofm/\Am_5\to C_0'':=C_0'/(t_1-1)=\GL_2.
\]
Since $\Am_6/\Am_5$ is (\latin{a priori} a quotient of)
the cyclic $\tGL_2$-module
$\tGL_2\ag_2$, the restriction of~$\partial_1$ maps it isomorphically onto
$(t_2-1)C_0''=\tGL_2$, proving Statement~\iref{ell.6}.

For the other statements, it suffices to estimate the number of generators.
With possible future applications in mind, we describe the structure of the
intermediate quotients in the form
$\text{(known module)}\onto\Am_k/\Am_{k-1}$. In fact, all these epimorphisms
are isomorphisms, see \autoref{rem.=} below.

In $\tAm\ofm/\Am_4$, one has
\[*
t_3=1,\qquad
\ag_1=\ag_3=0,\qquad
\ag'_2=-\cg'_3;
\]
the last relation follows from~\eqref{rel.a3}.
Thus, $\Am_5/\Am_4$ is generated by~$\cg'_3$, and
$\tAm\ofm/\Am_6$ is generated by $\cg_1$, $\cg_2$, $\cg_3$;
by~\eqref{rel.m} and \autoref{lem.b-module},
\[
\tGL_2\cg'_3\onto\Am_5/\Am_4,\qquad
\GL_1\cg_1\oplus\GL_2\cg_2\oplus\Z\cg_3\onto\tAm\ofm/\Am_6.
\label{epi.5,7}
\]
For the last summand $\Z\cg_3$, we use the fact that
\[*
(t_1-1)\cg_3=-(t_1-1)\ag_2=0\bmod\Am_6.
\]
Thus, $\ell(\tAm\ofm/\Am_6)\le2m+1$, and
Statements~\iref{ell.5} and~\iref{ell.7} are proved.

The module $\Am_3/\Am_1$ is generated by
$\bag_1$, $\bag_2$, $\bbg_1$, $\bbg_2$, $\bbg_3$,
and relations \eqref{rel.a1}, \eqref{rel.a2}, \eqref{rel.a3} imply
\[*
\bag_2=-\bbg_1,\qquad
\bag_1=-\bbg_2,\qquad
(t_1-1)(\bbg_3+\bag_2)=(t_2-1)\bag_1.
\]
We can retain three generators $\bbg_1$, $\bbg_2$, $\bbg_3$ only,
rewriting the last relation in the form
\[
(t_1-1)(\bbg_3-\bbg_1)+(t_2-1)\bbg_2=0.
\label{rel.A3}
\]
Note also that $\poly_m(t_3)\Am_3=0$, see~\eqref{rel.x}.

In $\Am_3/\Am_2$,
we have $(t_1-1)\bbg_3=(t_1-1)\bag_2=0$, hence also
$(t_1-1)\bbg_1=0$.
Then \eqref{rel.A3} implies $(t_2-1)\bbg_2=0$, and
\[
\tGL_3\bbg_1\oplus\tGL_3\bbg_2\oplus\tGL_3\bbg_3\onto\Am_3/\Am_2,
\label{epi.3}
\]
see \autoref{lem.b-module}.
This gives us Statement~\iref{ell.3}.

The
module $\Am_2/\Am_1$ is generated by~$\bbg'_1$
and~$\bbg'_3$.
By~\eqref{rel.x} and~\eqref{rel.bb}, we have
\[
\poly_m(t_i)(\Am_2/\Am_1)=0\quad\text{for all $i=1,2,3$}.
\label{rel.all}
\]
Relations~\eqref{rel.b2} and~\eqref{rel.A3} imply
$(t_1t_3-1)(\bbg'_3-\bbg'_1)=0$; using~\eqref{rel.bb}, this can be rewritten
as
$(t_3-t_2)\bbg'_3=(t_1-t_2)\bbg'_1$.
Let
\[*
\ug:=(t_3-t_2)\bbg'_3=(t_1-t_2)\bbg'_1
\]
and consider the cyclic
submodule $\Am_2'\subset\Am_2/\Am_1$ generated by~$\ug$.
By \autoref{lem.b-module},
\[
\tGL_2\bbg'_1\oplus\tGL_2\bbg'_3\onto(\Am_2/\Am_1)/\Am_2'.
\label{epi.2''}
\]
On the other hand, $\Am_2'\subset\GL\bbg'_1\cap\GL\bbg'_3$; hence,
$t_3\1=t_2=t_1\1$ on this module and, by \autoref{lem.b-module} again,
\[
\tGL_2\ug\onto\Am_2'\quad\text{if $m$ is odd}.
\label{epi.2'odd}
\]
This fact proves Statement~\iref{ell.2} in the case of $m$~odd.

If $m=2k$ is even, \eqref{epi.2'odd} still holds, but we need a stronger
statement.
Note that $\poly_m(t)$ is divisible by
$\poly_k(t^2)$. Furthermore,
one has a polynomial identity
\[
t^{m-2}\sum_{r=0}^{m-1}t^{1-r}\poly_r(t^2)=
 t\poly_{k-1}(t^2)\poly_m(t)+\poly_k(t^2),
\label{eq.poly}
\]
which is easily established by multiplying both sides by $t^2-1$.
On
the submodule $\Am'_2$ we have $t_2=t_1\1$, see~\eqref{rel.bb}; hence,
$s:=t_2t_1\1=t_2^2$. Then, representing $\ug$ in the form
$\ug=t_1(1-s)\bbg'_1$, we have
\[
t_2^{1-r}\poly_r(t_2^2)\ug=t_2^{1-r}t_1\poly_r(s)(1-s)\bbg'_1=
 t_1^r(1-s^r)\bbg'_1=(t_1^r-t_2^r)\bbg'_1,\quad r\in\Z.
\label{eq.ref}
\]
Summing up over $r=0,\ldots,m-1$ and using~\eqref{rel.all} and~\eqref{eq.poly}
at $t=t_2$, we conclude that
$\poly_k(t_2^2)\ug=0$, \ie,
\[
\GL_2\ug/\poly_k(t_2^2)\onto\Am_2'\quad\text{if $m=2k$ is even},
\label{epi.2'even}
\]
obtaining a stronger inequality
$\ell(\Am_2')\le\deg\poly_k(t^2)=m-2$.

The final inequality in the statement of the lemma is the sum of
items~\ref{ell.1}--\ref{ell.7}.
\endproof

\subsection{Proof of \autoref{th.main}}\label{proof.main}
We assume that $m\ge3$.
By \eqref{eq.section}, it suffices to show that
$\TT{\tln\ofm}=0$, where $\tln\ofm:=\ln\ofm+\rl_0\ofm$ is the divisor introduced
in \autoref{s.Fermat}. Since $\Fermat\ofm$ is simply connected, we
can use \autoref{th.reduction}, reducing the problem to proving the
inequality
$\ell(H_1(\Fermat\ofm\sminus\tln\ofm))\le\rank\KK{\tln\ofm}$.

According to~\cite{Shioda:Fermat1,Shioda:Fermat},
$\rank\sset_m=3(m-1)(m-2)+1+\Gd_m$.
On the other hand, $H_2(\tln\ofm)$ is the free abelian
group generated by the classes of the $3m^2$ lines and the additional class
$[\rl_0\ofm]$.
Hence, $\rank\KK{\tln\ofm}=9m-6-\Gd_m$, and the statement follows
from~\eqref{eq.H} and \autoref{lem.ell}.
\qed

\remark\label{rem.=}
It follows from the proof that all inequalities in the statement of
\autoref{lem.ell} are, in fact, equalities, \ie, no relation has been lost,
even though some relations were multiplied by non-units.
Furthermore, all epimorphisms
\eqref{epi.5,7},
\eqref{epi.3}, \eqref{epi.2''},
\eqref{epi.2'odd}, \eqref{epi.2'even} are isomorphisms.
\endremark

\remark\label{rem.<=}
We only use the inequality $\rank\sset_m\le3(m-1)(m-2)+1+\Gd_m$,
\ie, the fact that there is
\emph{at least} a certain number of relations between
the components.
In general, it would suffice to prove the inequality
$\ell(\Hm\of\Ga)\le\rank\KK{\dv\of\Ga}+\ls|G|-1$,
see \autoref{lem.H}.\mnote{ref added}
\endremark

\remark\label{rem.rank}
The rank\mnote{new remark} $\rank\sset_m$ can easily be computed directly,
by tensoring the module by~$\C$ and counting the irreducible summands, which
are all of dimension~$1$ (multi-eigenspaces of the three commuting
finite order operators $t_1$, $t_2$, $t_3$).
\endremark

\remark\label{rem.further}
By~\eqref{eq.T=Tors0},\mnote{new remark}
when computing the torsion, one can replace $\Hm\of\Ga$ with
the smaller module $\Bm\of\Ga$.
\latin{A posteriori}, $\Bm\ofm$ is the $\GL\ofm$-module
spanned by the three generators
$\cg_1,\cg_2,\cg_3$ subject to a single relation
\[*
(t_1-1)(t_3-1)\cg_1=(t_2-1)(t_3-1)\cg_2+(t_1-1)(t_3-1)\cg_3,
\]
see~\cite{degt:Delsarte}. In this form, some of the results of this paper
generalize to Fermat varieties of higher dimension, see~\cite{degt.Shimada}.
Note, though, that this one-relator presentation of $\Bm\of\Ga$ does not
extend to more general Delsarte surfaces; see~\cite{degt:Delsarte} for
further details.
\endremark

\subsection{A toy example}\label{s.toy}
In conclusion, we consider a very simple example, answering the generalized
question, see \autoref{problem}, in the special case of a covering ramified
over at most three lines.

\theorem\label{th.toy}
If the covering $\Fermat\of\Ga\to\Fermat$ is unramified over at least
one of the lines~$\rl_j$, $j=0,1,2,3$,
then $\TT{\dv\of\Ga}=0$.
\endtheorem

\proof
We can assume that the covering is unramified over~$\rl_3$, \ie, the
epimorphism $\Ga\:\GG\to G$ sends~$\g_3$ to~$0$.
Then, obviously, $\Ker\Ga=\Z\g_3\oplus(\GG_{12}\cap\Ker\Ga)$ and,
by \autoref{prop.H1}, we have $H_1(\Fermat\of\Ga)=0$, \ie,
\autoref{th.reduction} is applicable.

By~\eqref{rel.Ga}, we have $t_3=1$ on $\Hm\of\Ga$, and
relations~\eqref{rel.a1}, \eqref{rel.a2}, \eqref{rel.a3} become
\[*
(t_2-1)\ag_3=(t_1-1)\ag_3=0,\qquad
(t_1-1)(\cg_3+\ag_2)=(t_2-1)\ag_1.
\]
Introducing the generator $\ag_2':=\cg_3+\ag_2$ instead of~$\ag_2$, we
see that the submodule $\Bm\of\Ga\subset\Am\of\Ga$ introduced in
\autoref{s.finite} is a direct summand (as a $\GL$-module), and all relations
in $\Bm\of\Ga$ are $t_3=1$ and~\eqref{rel.b1}, \eqref{rel.b2}, \eqref{rel.b3}.
The three latter translate into independent relations
$(t_2-1)\cg_1=(t_1-1)\cg_2=(t_1t_2-1)\cg_3=0$, and $\Bm\of\Ga$ is a direct sum
of three group rings:
\[*
\Bm\of\Ga=\Z[G/\Ga(\g_2)]\cg_1\oplus\Z[G/\Ga(\g_1)]\cg_2\oplus
 \Z[G/\Ga(\g_1+\g_2)]\cg_3.
\]
By~\eqref{eq.T=Tors0},
one has $\TT{\dv\of\Ga}\cong\Tors\Bm\of\Ga=0$.
\endproof

\corollary[of~\eqref{eq.Shioda.Ga} and \autoref{th.toy}]\label{cor.toy}
If a covering $\pr\:\Fermat\of\Ga\to\Fermat$ as in \autoref{th.toy}
has degree~$m$ prime to~$6$, then $\Pic\Fermat\of\Ga=\SS{\dv\of\Ga}$.
\done
\endcorollary

\subsection{Proof of \autoref{th.pq} and \autoref{cor.pq}}\label{proof.pq}
\autoref{cor.pq}
is an immediate consequence of \autoref{th.pq} and the fact that
$\Pic\PQ_{p,q}$ is rationally generated by the classes of the lines,
see~\cite{Boissiere.Sarti}.
In view of \autoref{th.reduction}, the statement of \autoref{th.pq} is purely
homological, and we can deform $\PQ_{p,q}$ to the Fermat
surface~$\Fermat\ofm$; then, the $m^2$ lines in question deform to the
components of $\ln_1\ofm$, and $\sset_{p,q}=\SS{\ln_1\ofm}$.
Similar to~\eqref{eq.section}, the latter group equals $\SS{\tln_1\ofm}$,
where $\tln_1\ofm:=\ln_1\ofm+\rl_0\ofm$.

Patching $\ln_2\ofm$ and $\ln_3\ofm$, \cf.
\autoref{s.Fermat}, we conclude that
\[*
\TT{\tln_1\ofm}=\Tors\tAm'\ofm,\qquad
\tAm'\ofm:=\tAm\ofm/(\GL\cg_2+\GL\cg_3).
\]
Filtering this module as in \autoref{lem.ell} and analyzing the proof of the
lemma, we see that Statements~\iref{ell.1}, \iref{ell.4}, and~\iref{ell.6}
hold without change, whereas the other statements can be rewritten as
follows:
\roster
\addtocounter{enumi}1
\item\label{ell'.2}
$\ell(\Am_2/\Am_1)=0$ due to \eqref{rel.A3},
\item\label{ell'.3}
$\ell(\Am_3/\Am_2)\le m-1$, see~\eqref{epi.3},
\addtocounter{enumi}1
\item\label{ell'.5}
$\ell(\Am_5/\Am_4)=0$, see~\eqref{epi.5,7},
\addtocounter{enumi}1
\item\label{ell'.7}
$\ell(\Am_7/\Am_6)\le m$, see~\eqref{epi.5,7}.
\endroster
Summing this up, we obtain $\ell(\tAm'\ofm)\le m^3+2m-2$.
On the other hand,
one has $\rank\SS{\tln_1\ofm}=(m-1)^2+1$, see~\cite{Boissiere.Sarti};
hence, $\rank\KK{\tln_1\ofm}=2m-1$ and, as in \autoref{proof.main}, we
conclude that $\tAm'\ofm$ is a free abelian group.
\qed

\let\.\DOTaccent
\def\cprime{$'$}
\bibliographystyle{amsplain}
\bibliography{degt}


\end{document}